\documentclass[3p,times,procedia,number]{elsarticle}
\usepackage{ecrc}
\graphicspath{ {./img/} }

\volume{00}

\firstpage{1}

\journalname{Procedia Engineering}

\runauth{N. Barral et al.}

\jid{proeng}

\usepackage{amssymb}
\usepackage[utf8]{inputenc}
\usepackage{stmaryrd,latexsym,multirow,array}
\usepackage{float}
\usepackage{subfig}
\usepackage{amssymb,amsbsy,amsmath,amsfonts,epsfig}
\usepackage{amsthm}
\usepackage{graphicx}
\usepackage{algorithm}
\usepackage{url}
\usepackage{colortbl}

\definecolor{blue}{RGB}{0,0,255}

 \biboptions{sort&compress}

\usepackage[figuresright]{rotating}

\usepackage{biblioSty}

\begin{document}

\begin{frontmatter}



\dochead{25th International Meshing Roundtable}

\title{Anisotropic mesh adaptation in Firedrake with PETSc DMPlex}


\author[a]{Nicolas Barral}
\author[b]{Matthew G. Knepley}
\author[a]{Michael Lange}
\author[a]{Matthew D. Piggott}
\author[a]{Gerard J. Gorman\corref{cor1}}

\address[a]{Department of Earth Science and Engineering, Imperial College, South Kensington Campus, London SW7 2AZ, UK}
\address[b]{Computational and Applied Mathematics,
Rice University,
6100 Main MS-134,
Houston, TX 77005,
USA}

\begin{abstract}

Despite decades of research in this area, mesh adaptation capabilities are 
still rarely found in numerical simulation software. We postulate that the 
primary reason for this is lack of usability. Integrating mesh adaptation 
into existing software is difficult as non-trivial operators, such as error 
metrics and interpolation operators, are required, and integrating available 
adaptive remeshers is not straightforward. Our approach presented here is to 
first integrate Pragmatic, an anisotropic mesh adaptation library, into 
DMPlex, a PETSc object that manages unstructured meshes and their 
interactions with PETSc's solvers and I/O routines. As PETSc is already 
widely used, this will make anisotropic mesh adaptation available to a much 
larger community. As a demonstration of this we describe the integration of 
anisotropic mesh adaptation into Firedrake, an automated Finite Element based
system for the portable solution of partial differential equations which 
already uses PETSc solvers and I/O via DMPlex. We present a proof of concept 
of this integration with a three-dimensional advection test case.

\end{abstract}

\begin{keyword}
Anisotropic Mesh Adaptation ; Finite Element ; PETSc ; Global Fixed-Point Adaptation Algorithm 




\end{keyword}
\cortext[cor1]{Corresponding author. Tel.: +44 (0)20 7594 9985.}
\end{frontmatter}

\email{g.gorman@imperial.ac.uk}





\section{Introduction}
 
As the size of the numerical simulations required by both the industry and 
research sectors continues to grow, anisotropic mesh adaptation is an 
efficient means to reduce the CPU time of the computations while improving 
their accuracy~\cite{Pain-2001,Compere-2008,Alauzet-2016}. However, it is 
still not widely adopted in the solvers community due to the difficulty of 
integrating it to partial differential equation solver frameworks. Various 
adaptive mesh processes can be found in the literature, but they all remain 
complex to set up. They require several complex steps: computation of an 
element size and orientation map, generation of a mesh with respect to this 
map and solution transfers, which are often implemented in different codes 
that have to be combined. Adaptation for transient problems is even more 
complex, as these steps often have to be repeated many times.

In PETSc, the widely used scientific library providing data structures and 
routines for the parallel solution of partial differential 
problems~\cite{petsc-user-ref}, the mesh management library 
DMPlex~\cite{Knepley-2009} aims at making
application code developers' lives easier by providing them with a wide set 
of tools to manipulate unstructured meshes, thus avoiding the need to write 
such specialized routines. 
Integrating Pragmatic~\cite{Gorman-2015}, an anisotropic mesh adaptation 
library, into DMPlex is a first step towards making mesh adaptation available 
to a larger community, as it offers the possibility to easily generate and use 
adapted meshes within application codes which make use of PETSc.
As an example, Firedrake~\cite{Rathgeber-2015}, a system for the portable 
automated solution of partial differential equation based problems, is one 
such code considered here. Providing routines for mesh adaptation in this 
framework would make it even more accessible.

The purpose of this research note is to present the integration of Pragmatic 
to Firedrake via DMPlex and the implementation of a full mesh adaptation 
algorithm with Firedrake, which are the first steps towards fully automated 
mesh adaptation. We first present the three codes
involved, and show how they are integrated. The chosen adaptation strategy is 
then described, and demonstrated using a 3D advection problem.

\section{Presentation of the codes involved }


\paragraph{Firedrake}

Firedrake is a code for the automated solution of Finite Element
problems~\cite{Rathgeber-2015}. It is based on a strict division into different
abstraction layers, each layer being concerned with a specific task: definition
of the problem, local discretization defining the data structures and kernels
to compute the solution, and the parallel execution of these kernels.  Different
kinds of optimization can be applied at each level, from caching of
mathematical forms to compiler-level optimization in the kernels, including
data sorting. The definition of the problem is done using a domain-specific
language (DSL) for the specification of partial differential equations in
variational form. This allows the user to define very simply their problem, with
a language close to the underlying mathematical description. A wide range of finite element
types and degrees are supported. Efficient parallel execution of the numerical
kernels on unstructured mesh data is performed by a framework where parallel
loops are expressed at a high level, and which uses its own
representation of the data as sets and maps. Finally, the resulting
linear or non-linear systems are solved via the PETSc library, which is
also used for the handling of unstructured meshes via DMPlex.

\paragraph{PETSc and DMPlex}

\begin{figure}[h!]
  \centering
    \includegraphics[width=0.52\linewidth]{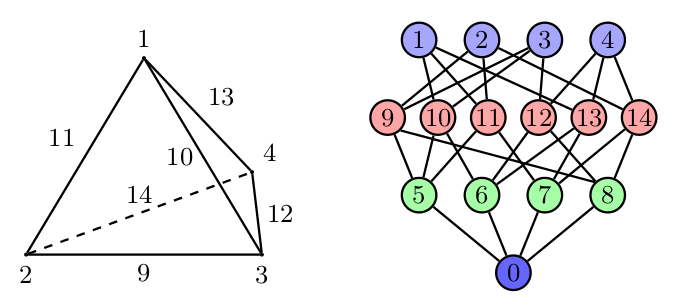} 
    \caption{\label{fig-dmplex} Example of DMPlex representation for a single tetrahedron. }
\end{figure}

DMPlex is the PETSc object which manages unstructured meshes using a Hasse
diagram representation and directed acyclic graph data
structure~\cite{Knepley-2009}. The library has a wide variety of tools for mesh
traversal and manipulation, including submesh extraction, partitioning and
distribution, I/O, and parallel load balancing. DMPlex uses a uniform
topological interface, meaning that vertices, edges, facets and cells are
treated equally as points of the graph representation (see Fig.~\ref{fig-dmplex}). 
This allows algorithms
which operate independent of dimension, cell type, cell mixture, or
partitioning. A key point is that topological operations can be directly
translated to graph operations, which is ideal for optimization and parallelization.
DMPlex also provides an abstraction for associating function spaces with pieces
of the mesh, allowing the user to construct approximations like finite elements
and finite volumes which are made of small spaces pieced together. Firedrake
uses DMPlex for mesh representation, to organize traversals in the global
approximation space, and to interface with advanced solvers, such as
unstructured multigrid and block preconditioners. The way DMPlex enables
efficient handling of unstructured meshes is detailed in~\cite{Lange-2015}.

\paragraph{Pragmatic}\label{sec-pragmatic}

Pragmatic is a 2D and 3D anisotropic local remesher, that generates a unit 
mesh with respect to a prescribed metric field~\cite{Gorman-2015}. The 
adapted mesh is obtained from the input mesh through a series of local mesh 
manipulations: iterative applications of refinement, coarsening and edge/face
swapping optimise the resolution and the quality of the mesh, and 
quality-constrained Laplacian smoothing fine-tunes the mesh quality. An 
interface between DMPlex and Pragmatic was developed in this work. This 
enables PETSc users to adapt an unstructured mesh to a specified metric field
with only one call. The DMPlex object is converted into Pragmatic data 
structures, and a new DMPlex is created from its output. The advantage of 
this approach is that it allows application codes to interface with mesh 
adaptation through an I/O like interface, thereby making the integration of 
mesh adaptation much less intrusive in the application code.

\section{Transient adaptation algorithm }
\label{sec-ptfxA}

To demonstrate the solution we consider the global fixed-point adaptation
algorithm described in~\cite{Barral-thesis}
, based on: 
\begin{itemize}
  \item the subdivision of the simulation interval into $n_{adap}$
sub-intervals: $[0,\,T]~=~ [0=t^1,\,t^2]\cup \ldots \cup[t^i,\,t^{i+1}]\cup
\ldots \cup[t^{n_{adap}},\,t^{n_{adap}+1} = T]$, on each of which we consider
a unique adapted mesh,
  \item an iterative process to converge the mesh/solution couple and to
address prediction of the solution issues. 
\end{itemize}
An overview of the algorithm is presented in Algorithm~\ref{alg-uaa}. On each
sub-interval $[t^i, t^{i+1}]$, we consider a mean Hessian of the solution $u$ :
${\mathbf H}^i_{u}({\bf x}) = \int_{t^{i}}^{t^{i+1}} |{H}_u({\bf x},t)| \,
\mathrm{d}t$. The $d$-dimensional adapted mesh on sub-interval $i$ is generated
from the following metric, which has been proven to be optimal for the control of
the $L^p$ interpolation error:
\begin{equation}\label{eq-Metric_MLp_FixedTau_SubInt}
{\mathcal M}^i_{L^p}(\mathbf{x})   =   {\mathcal N}_{st}^{\frac{2}{d}} 
\biggl( \sum_{j=1}^{n_{adap}} {\mathcal K}^j \, \Bigl( \int_{t^{j}}^{t^{j+1}} \tau(t)^{-1} \mathrm{d}t \Bigr)^{\frac{2p}{2p+d}}  \biggr)^{-\frac{2}{d}} 
\Bigl( \int_{t^{i}}^{t^{i+1}} \tau(t)^{-1} \mathrm{d}t \Bigr)^{-\frac{2}{2p+d}}  (\det {\mathbf H}_u^i({\bf x}))^{-\frac{1}{2p+d}} \, {\mathbf H}_u^i({\bf x})\,,  
\end{equation}
where $\displaystyle {\mathcal K}^i = \int_{\Omega} \left( \det {\mathbf
H}_u^i({\bf x}) \right)^{\frac{p}{2p+3}} \mathrm{d} \mathbf{x} $, $\tau(t)$ is
a function of time specifying the time-step and $\mathcal{N}_{st}$ is the
target space-time complexity (the number of space-time vertices). 
${\mathcal N}_{st}^{\frac{2}{d}} \biggl(
\sum_{j=1}^{n_{adap}} {\mathcal K}^j \, \Bigl( \int_{t^{j}}^{t^{j+1}}
\tau(t)^{-1} \mathrm{d}t \Bigr)^{\frac{2p}{2p+d}}  \biggr)^{-\frac{2}{d}} $ 
is a global normalization term that
enables an optimal distribution of the number of vertices on the 
sub-intervals. It requires that the problem is solved on the whole simulation 
interval before the metrics for each sub-interval can be computed. The steps 
of Algorithm~\ref{alg-uaa} are detailed in the following Section.

\begin{algorithm}[!ht]
\caption{Mesh Adaptation Algorithm}
Given: Initial mesh and solution, target space-time complexity \\
For $j=1,n_{ptfx}$ \quad \quad \quad \quad  {\em // Fixed-point loop to converge the mesh/solution couple problem}\\
\vspace*{-4mm}
\begin{enumerate}
  \vspace*{-4mm}
  \item For $i=1,n_{adap}$ \quad \quad  {\em // Adaptive loop to advance the solution in time on time frame $[0,T]$}
  \vspace*{-4mm} 
  \begin{enumerate}
    \item Interpolate the solution from the previous sub-interval;
    \item Compute the solution on sub-interval;
    \item Compute sub-interval averaged Hessian from Hessian samples;
  \end{enumerate}
  \vspace*{-4mm}
  EndFor
  \item Compute global normalization term and all sub-interval metrics;
  \item Generate all sub-interval adapted meshes;
\end{enumerate}
\vspace*{-4mm}
EndFor
\label{alg-uaa}
\end{algorithm}

\section{Implementation of the algorithm}

\paragraph*{Metric computation}

The metrics from Equation~(\ref{eq-Metric_MLp_FixedTau_SubInt}) are based on
the mean Hessian of the solution. Instead of
computing the Hessian for each solver time-step, we only average samples,
typically 20 per sub-interval. The Hessian, $H$, is computed using a classic
Galerkin method. We desire $H$ such that: $H = D^2 u$, which can be written in weak
form, for a test function $\tau$ in the appropriate function space $\Sigma$:
\begin{equation}
    \int_\Omega \tau:(\sigma + D^2u)\,\mathrm{d}x - \int_{\partial\Omega} \tau:(\mathbf{n} \otimes \nabla u)\,\mathrm{d}s = 0, \quad \forall \tau \in \Sigma \,.
\end{equation}
This can be easily written in the high-level language of Firedrake, and then 
is automatically solved.
The samples are then averaged on the fly, and the resulting mean
Hessian for the sub-interval is stored. At the end of the simulation-interval,
the normalization term can be computed and the final metrics for each
sub-interval obtained.

\paragraph*{Generation of adapted meshes}

The meshes are generated by Pragmatic, whose integration with Firedrake via
DMPlex is explained in Section~\ref{sec-pragmatic}. For each sub-interval, the
DMPlex object corresponding to the current mesh is sent to PETSc together with
the metric field, and PETSc returns a DMPlex object corresponding to the new
adapted mesh.

\paragraph*{Solution transfer}

Between two sub-intervals, the solution has to be transferred from one mesh to
the following. For now, we use the point evaluation mechanism from Firedrake:
given any point in the domain, it provides the value of a given function at
this point. Since we are considering Lagrange $P^1$ finite elements, a linear
interpolation is performed, in the computational space.

\section{Numerical example}

We consider the 3D advection of a bubble, first introduced
in~\cite{Leveque-1996}. This case is usually used as a demonstrator for
interface-tracking methods. Here, we show that mesh adaptation allows us to
limit considerably the spurious diffusion of the solution.

We consider the scalar advection equation for a quantity $u$ in a 3D domain $\Omega$:
\begin{equation}
 \label{eq-adv}
 \frac{\partial u}{\partial t} + \nabla\cdot \left(\mathbf{v} u \right) = 0\,,
\end{equation}
where $\mathbf{v}$ is a velocity field defined on $\Omega$. The computational
domain is a $[0,1]\times[0,1]\times[0,1]$ cube. At $t=0$, $u=0$ everywhere 
except in a ball of radius $0.35$ centered in $(0.35,0.35,0.35)$, which models 
the bubble and where $u=1$. The following velocity field is considered, so the 
bubble quickly becomes distorted:
\begin{equation}
  \label{eq-vel}
  \mathbf{v}(x,y,z,t) = \left\{ \begin{aligned}
    2\sin^2(\pi x) \sin(2\pi y) \sin(2\pi z) \cos(2\pi t/T) \\
    -\sin(2\pi x) \sin^2(\pi y) \sin(2\pi z) \cos(2\pi t/T)  \\
    -\sin(2\pi x) \sin(2 \pi y) \sin^2(\pi z) \cos(2\pi t/T)
  \end{aligned}\right.\,,
\end{equation}
where $T$ is the period. Here $T=6$, and the simulation is run until $t=1.5$.


Equations~(\ref{eq-adv})~and~(\ref{eq-vel}) are solved using a Lagrange $P^1$ 
Finite Element method. A SUPG stabilization method is used, with a 
Crank-Nicolson scheme employed for advancing in time. This is easily 
achieved thanks to Firedrake: one just has to write the variational 
formulation in the high-level language, and the corresponding low-level solver 
is automatically generated. 

\begin{figure}[ht!]
  \centering
  \begin{tabular}{ccccc}
    \includegraphics[width=0.174\linewidth]{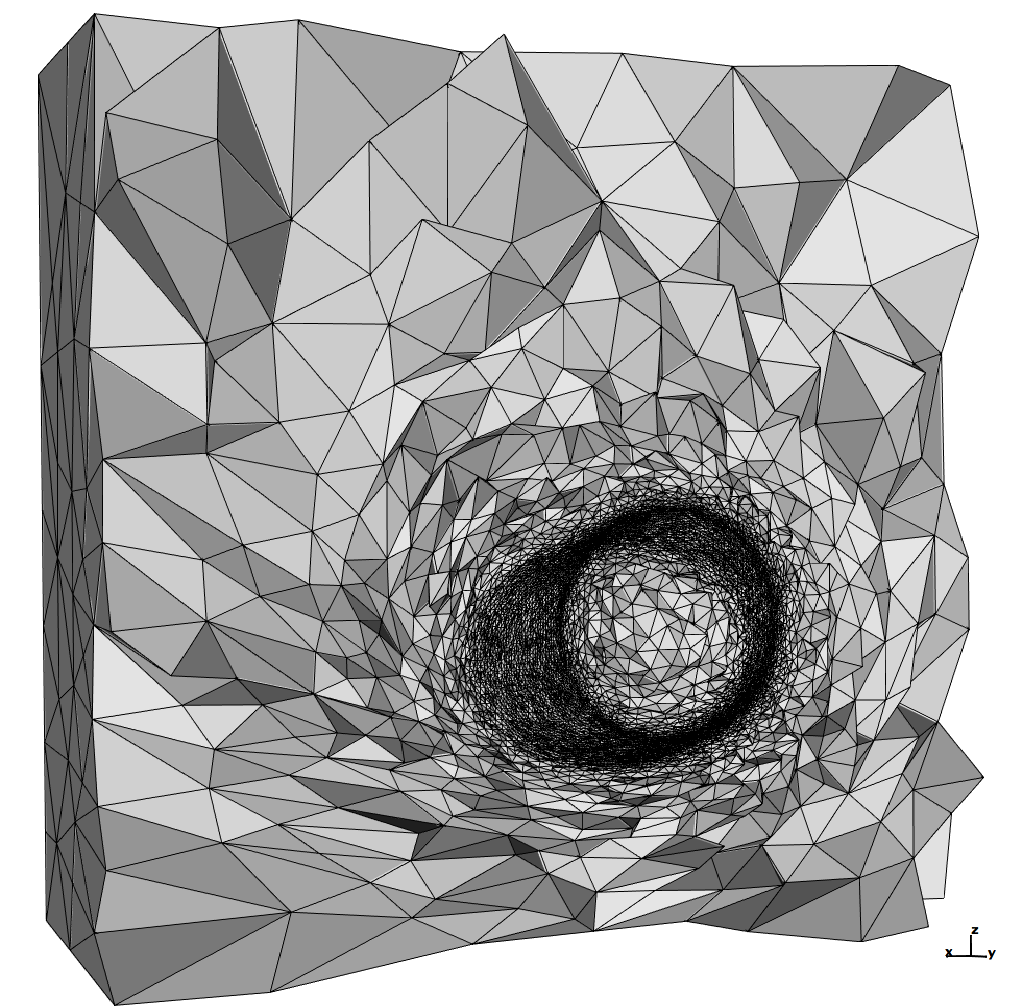} &
    \includegraphics[width=0.174\linewidth]{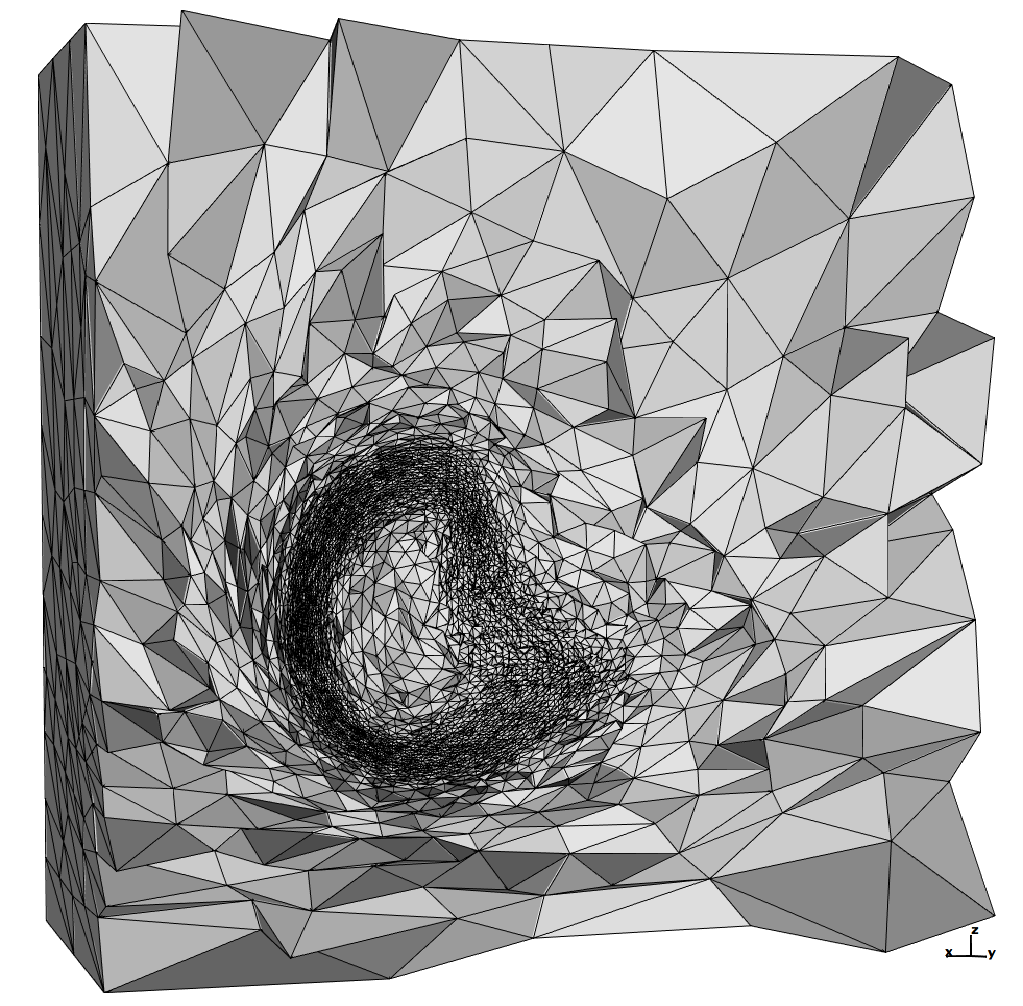} &
    \includegraphics[width=0.174\linewidth]{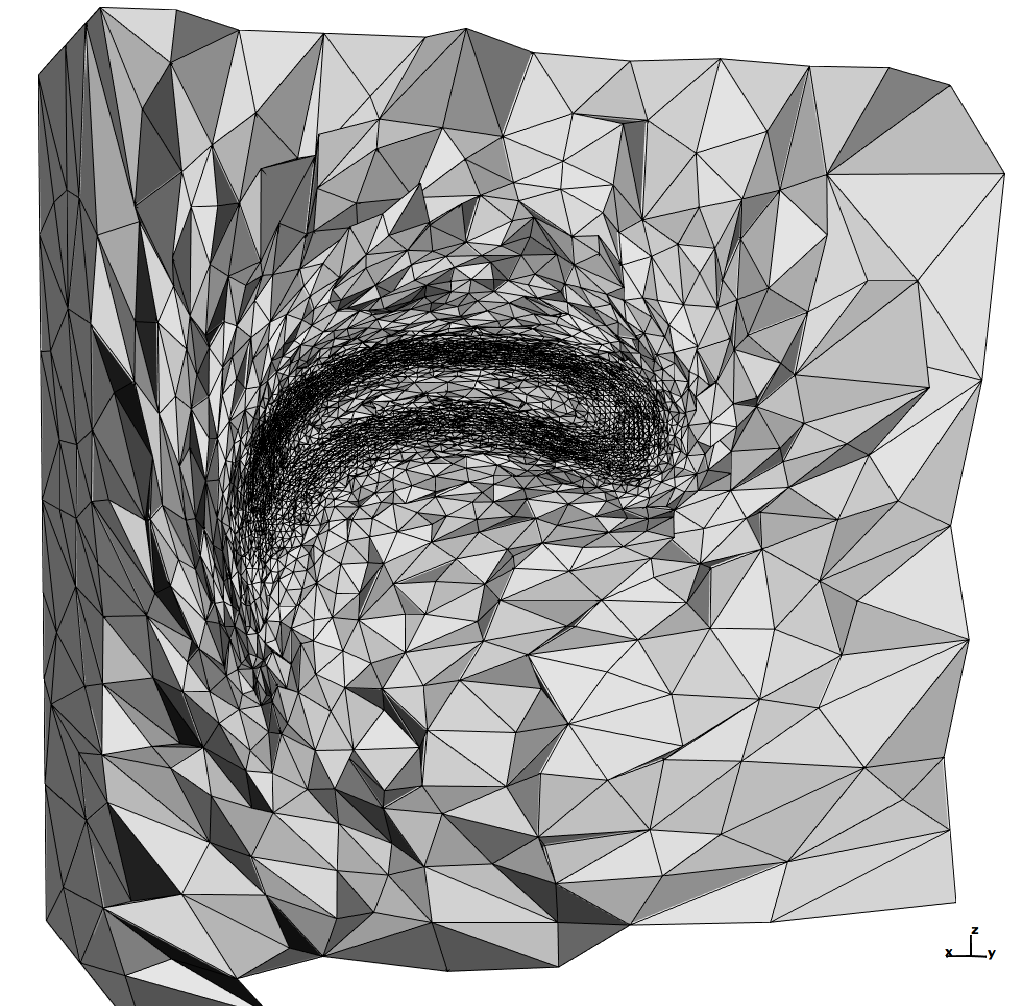} &
    \includegraphics[width=0.174\linewidth]{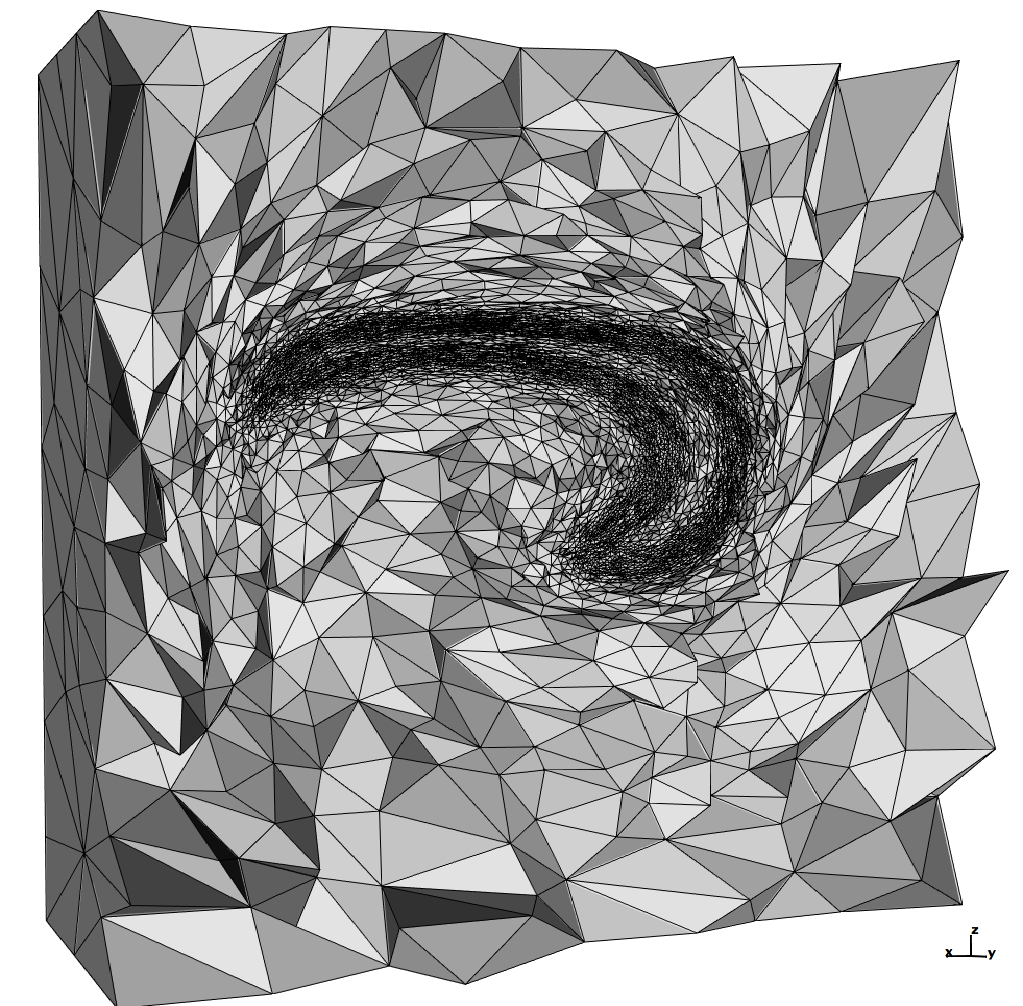} &
    \includegraphics[width=0.174\linewidth]{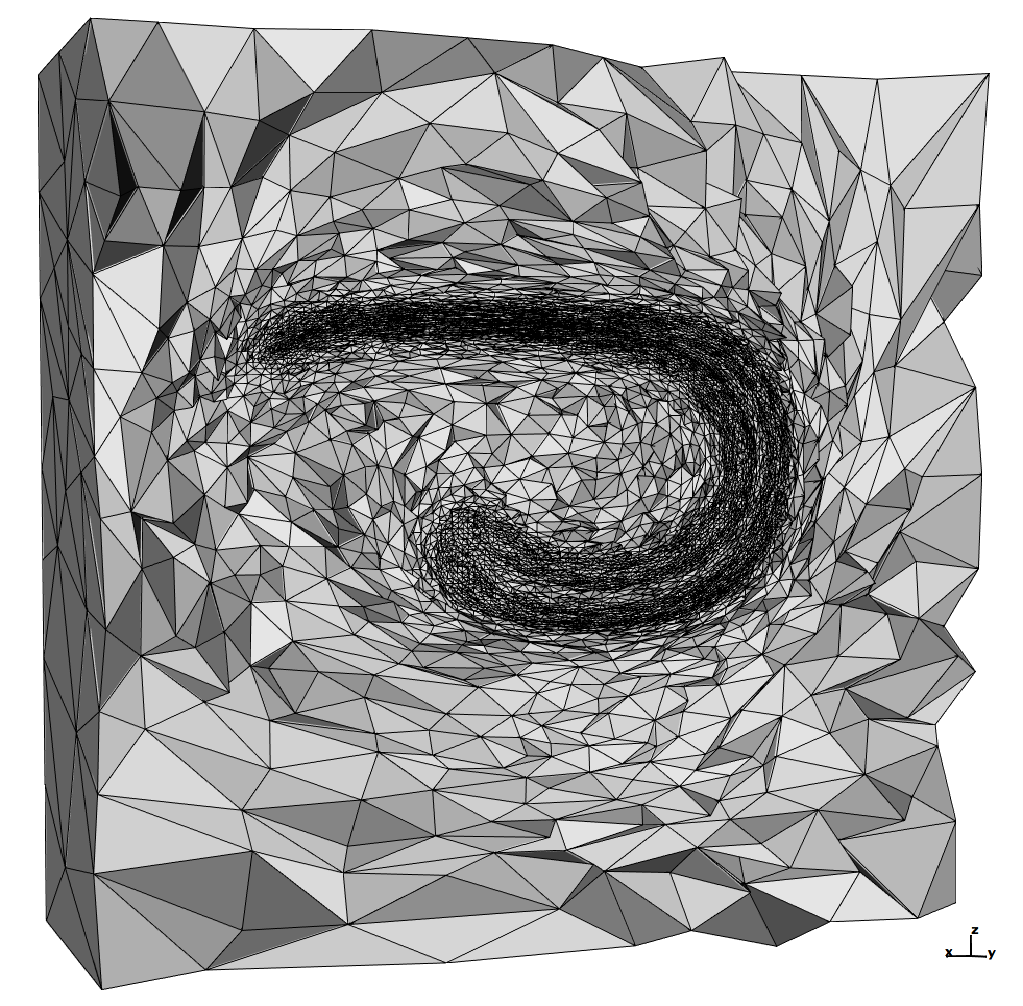} \\
    \includegraphics[width=0.174\linewidth]{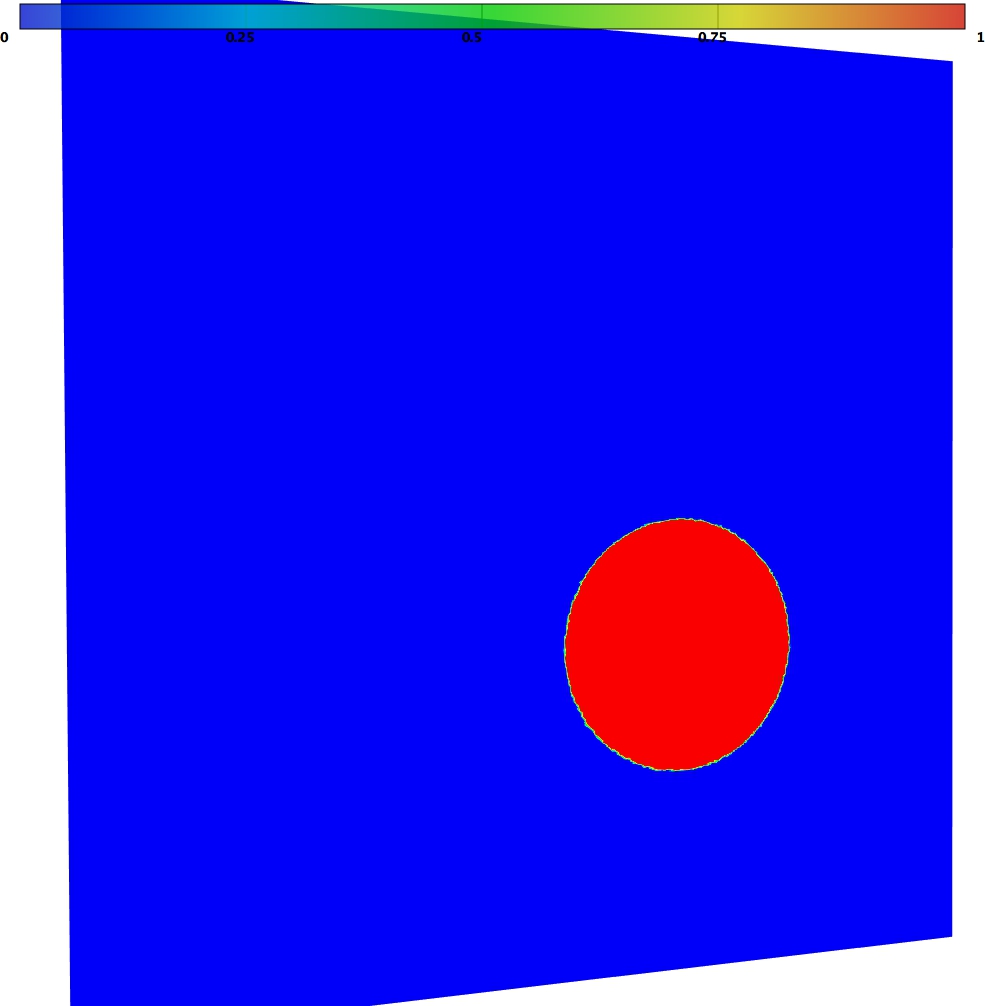} &
    \includegraphics[width=0.174\linewidth]{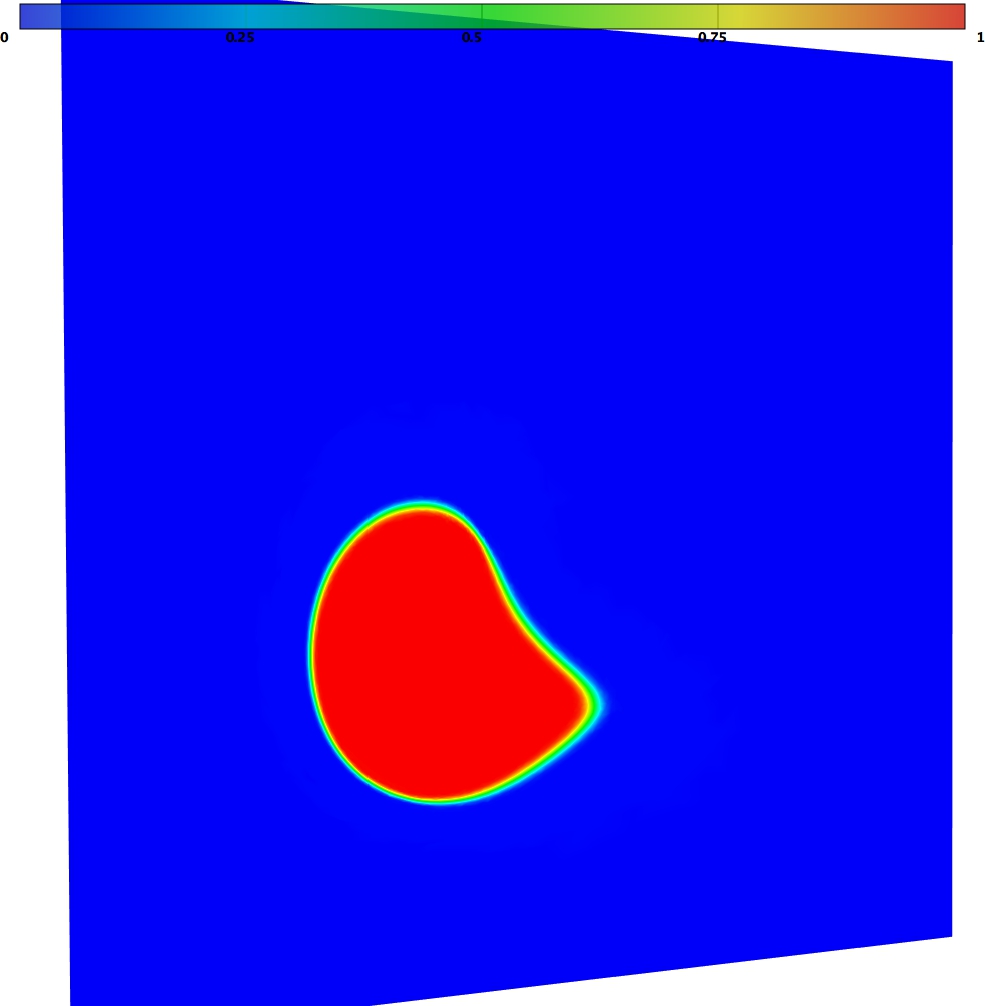} &
    \includegraphics[width=0.174\linewidth]{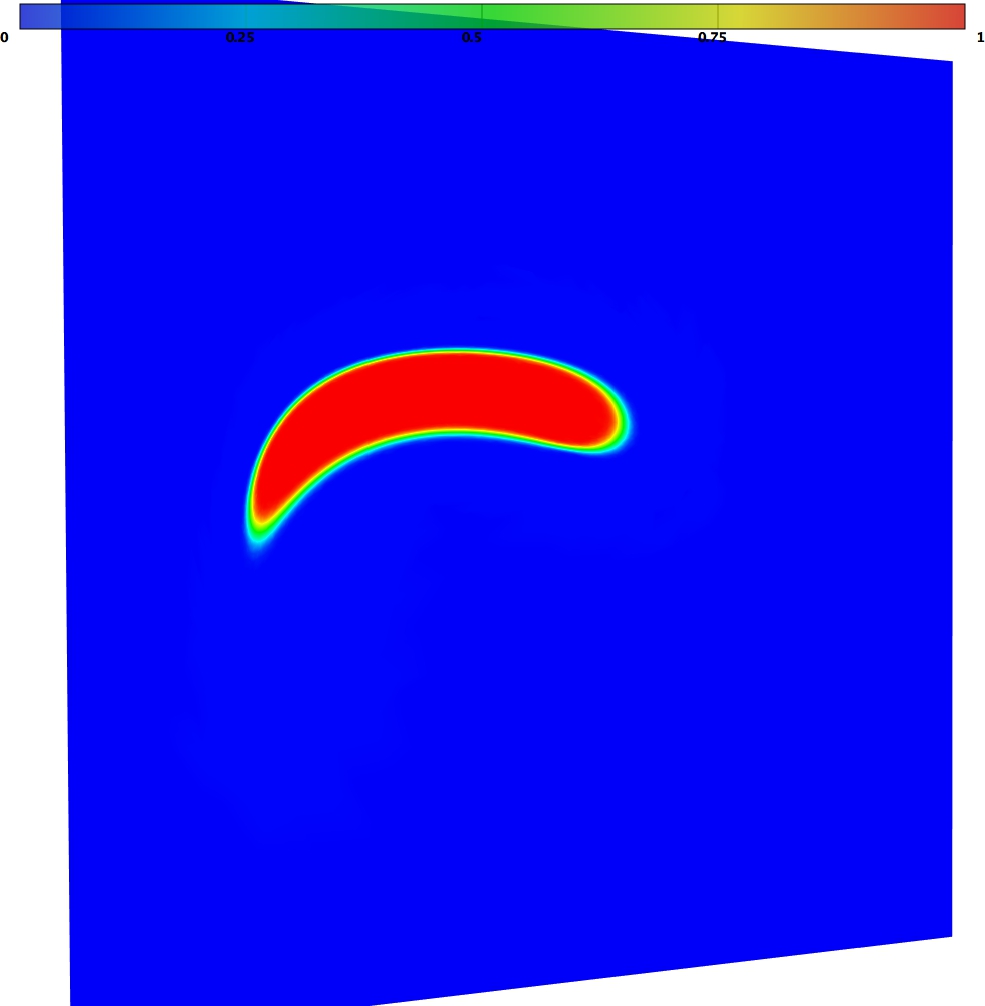} &
    \includegraphics[width=0.174\linewidth]{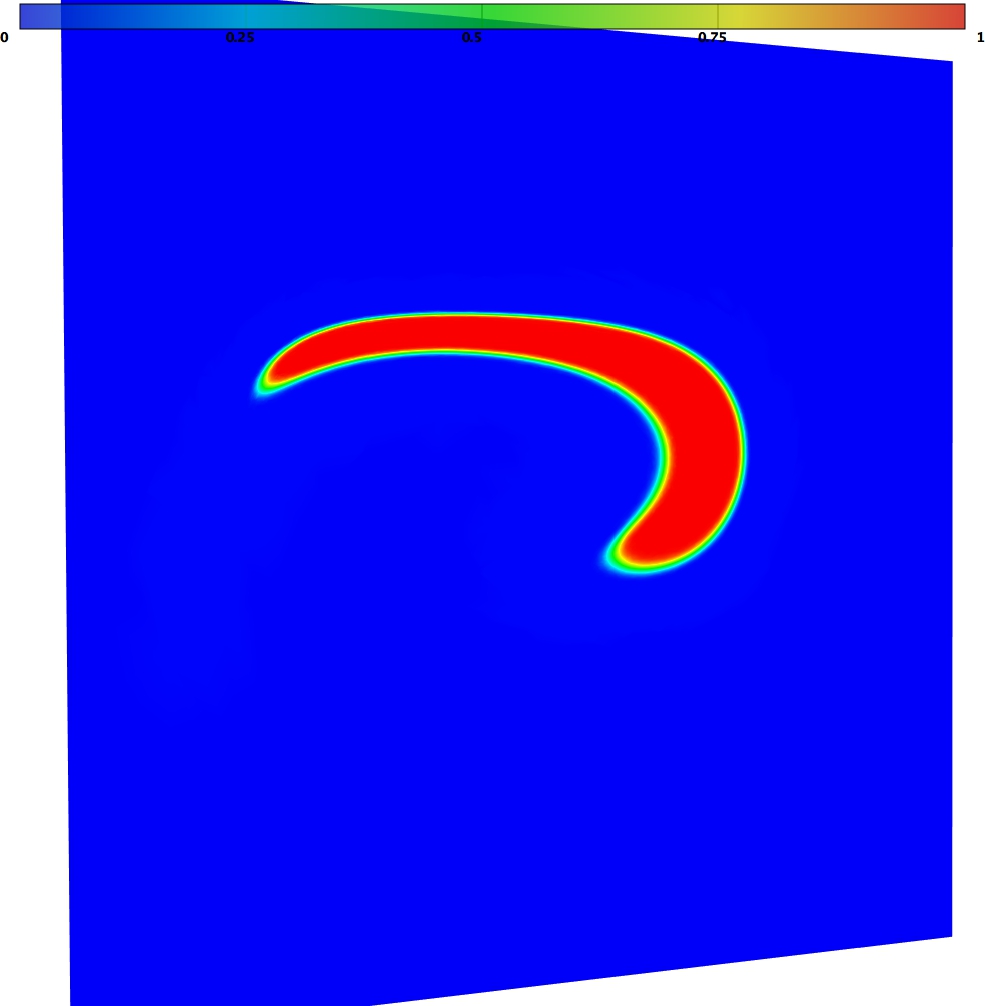} &
    \includegraphics[width=0.174\linewidth]{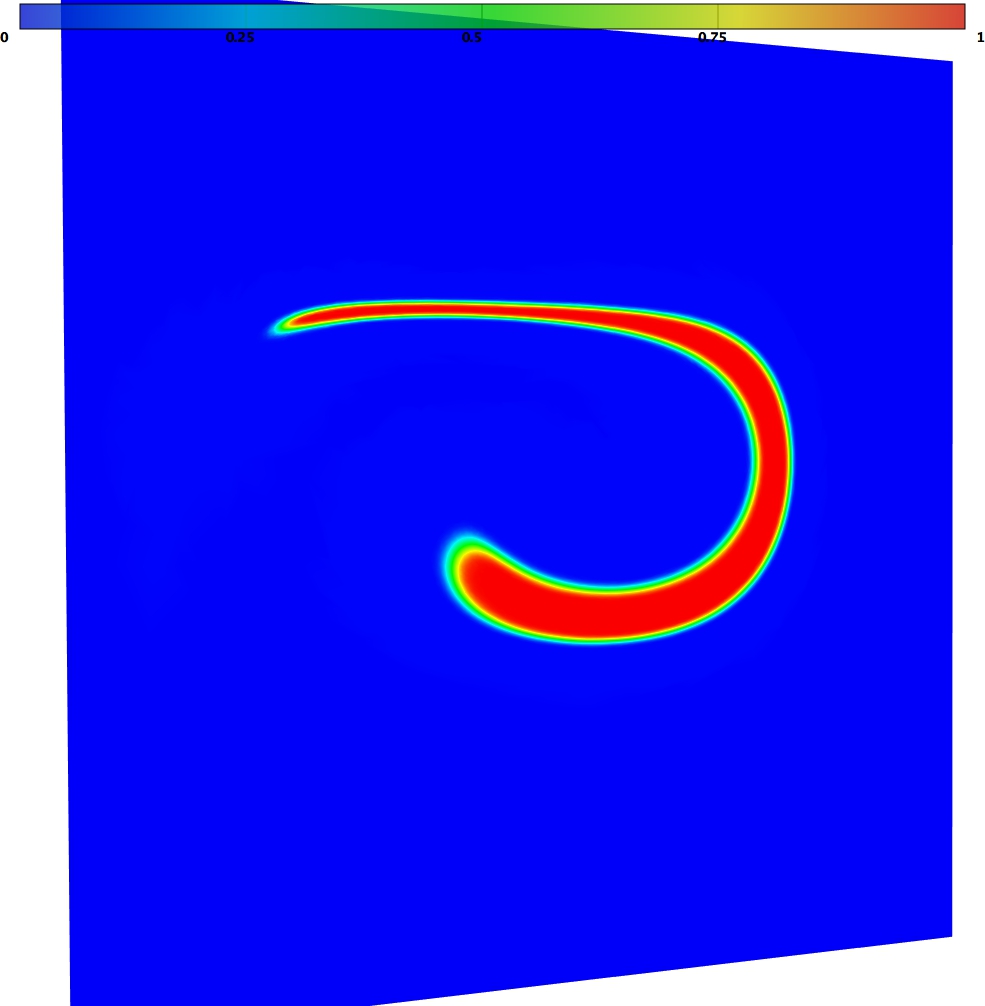} \\
    \includegraphics[width=0.174\linewidth]{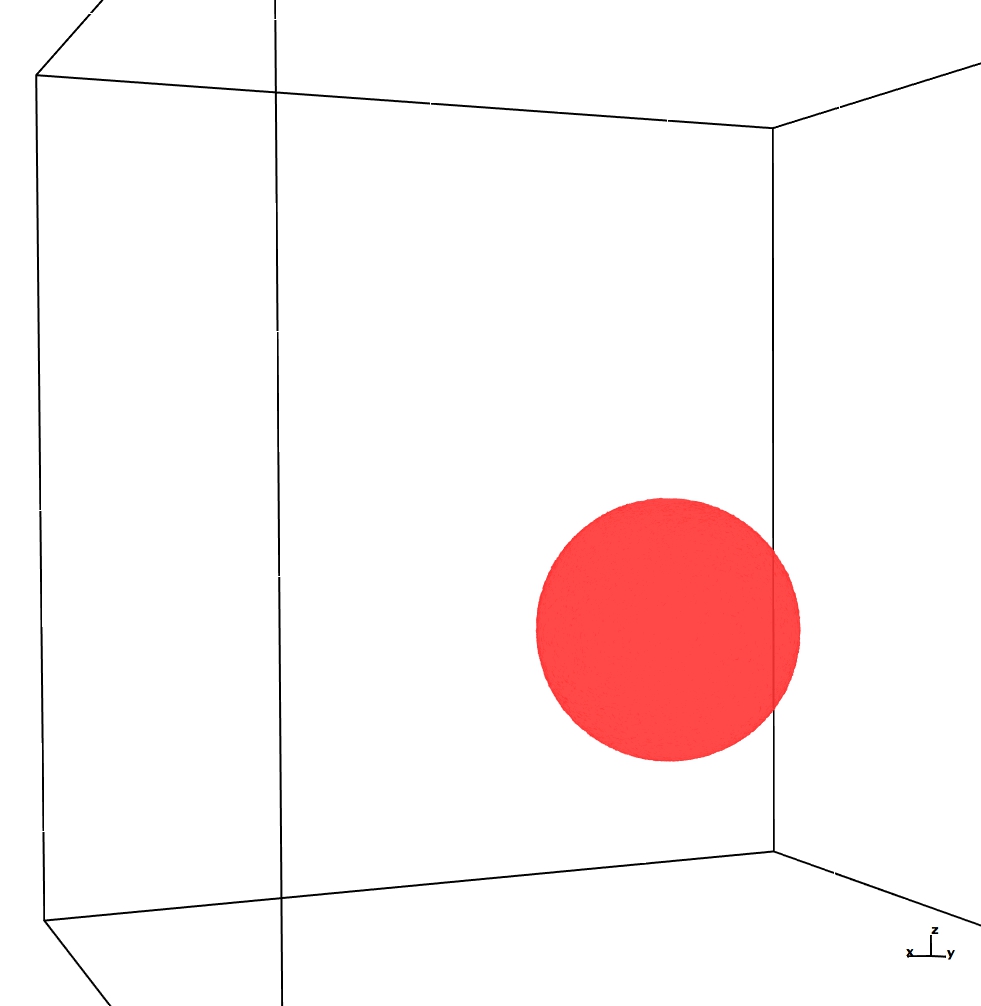} &
    \includegraphics[width=0.174\linewidth]{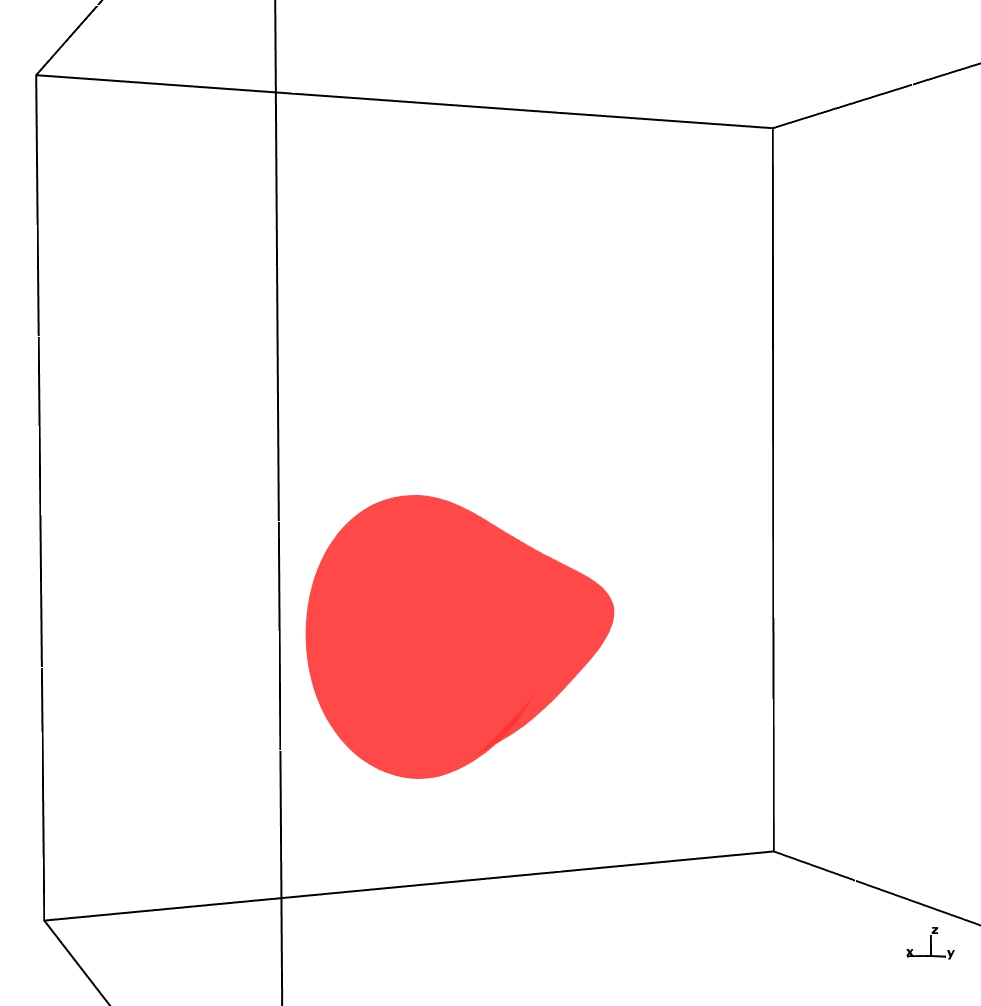} &
    \includegraphics[width=0.174\linewidth]{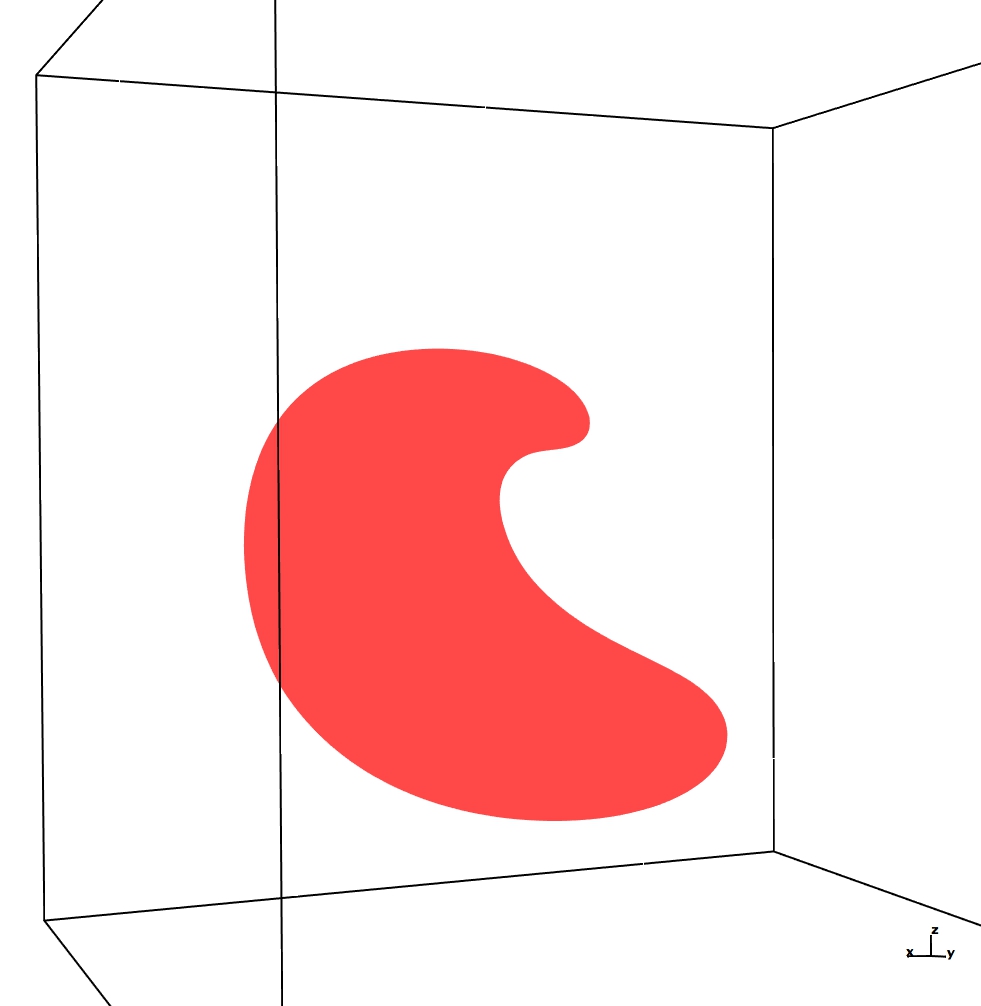} &
    \includegraphics[width=0.174\linewidth]{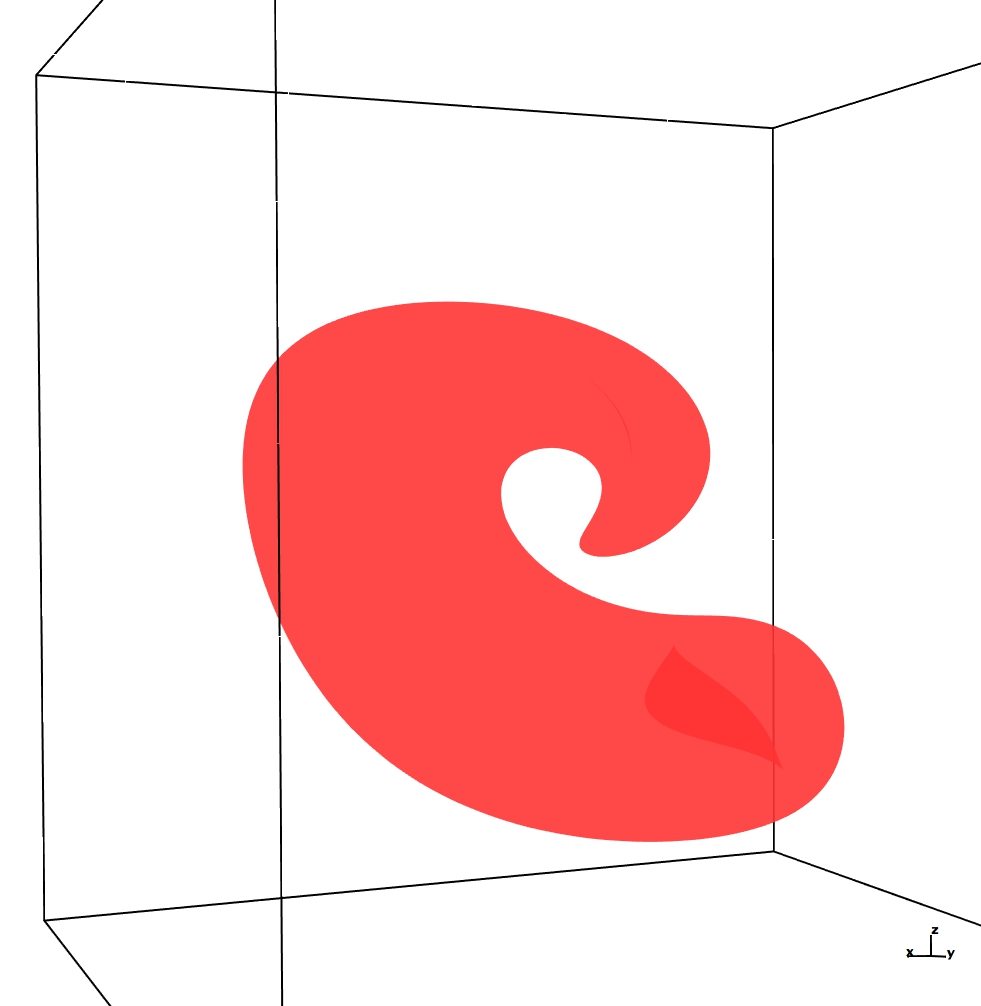} &
    \includegraphics[width=0.174\linewidth]{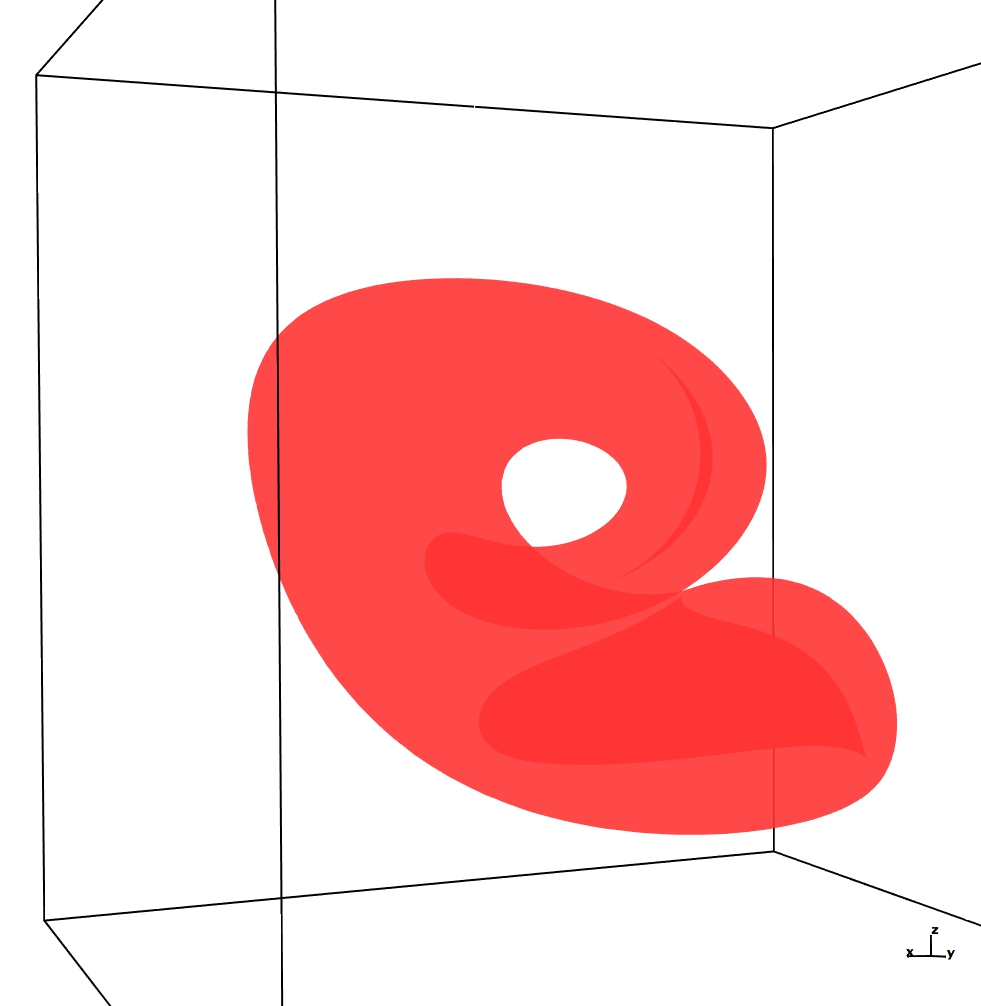} \\
  \end{tabular}
    \caption{\label{fig-bubble3d} Cuts in the mesh (top) and solution (middle)
     at dimensionless times $t=0,0.24, 0.54, 0.84\ \mathrm{ and }\ 1.5$. 
     Isosurfaces of the solution (bottom) at the same times are also given.}
\end{figure}

The global fixed-point algorithm from Section~\ref{sec-ptfxA} was used, with
$n_{adap} = 25$ sub-intervals and 3 global iterations. This way, only
25 adapted meshes are used for the whole simulation interval, and only
50 adapted meshes are generated. Snapshots of the meshes and solutions at the 
last fixed-point iteration are shown in figure~\ref{fig-bubble3d}. The meshes 
have an average size of $406,571$ vertices, ranging from $273,061$ to $557,586$. One 
can clearly see the adapted band-shaped regions, as well as the accuracy of 
the solution, that exhibits little numerical diffusion. 

The whole simulation was run serially, as the code is not yet fully 
parallelized. However, it was run in 32 hours on one core, which shows the 
efficiency of the approach and the potential power when run with more cores. 

\section{Conclusion}

We have shown how a transient mesh adaption algorithm was successfully 
implemented using Firedrake, PETSc DMPlex and Pragmatic. The long-term 
goal is to fully integrate mesh adaptivity within Firedrake's core, so that it 
can be used by application developers with only a few lines of code. The main 
ongoing improvement tracks include the parallelization of the whole process 
and the extension of this work to more complex finite element spaces. Most 
tasks are already automatically parallelized in Firedrake, what remains is 
the parallelization of the
interface between DMPlex and Pragmatic, for which parallel I/O support will 
have to be added in DMPlex. Extending adaptivity to higher degree or other 
kinds of Finite Elements notably requires improving the solution transfer 
procedure, and we aim at using the Supermesh library~\cite{Farrell-2009} to 
this end.

\section*{Acknowledgements}
This work was supported by EPSRC grant EP/L000407/1.




\bibliographystyle{plain}
{\small
{\setlength{\baselineskip}{0.7\baselineskip}
\bibliography{biblio}
}
}


\end{document}